\documentclass{article}
\usepackage{amssymb}  
\usepackage{spconf,amsmath,graphicx,hyperref}
\usepackage{enumitem}
\usepackage{subcaption} 

\title{Fractional harmonic transform on point cloud manifolds}
\name{Jiamian Li$^{1}$, Bing-Zhao Li$^{1*}$}
\address{$^{1}$School of Mathematics and Statistics, Beijing Institute of Technology, Beijing, China\\
E-mail: krown0205@163.com, li\_bingzhao@bit.edu.cn}

\begin{document}
%
\maketitle
\begin{abstract}
Three-dimensional point clouds can be viewed as discrete samples of smooth manifolds, allowing spectral analysis using the Laplace-Beltrami operator (LBO). However, the traditional point cloud manifold harmonic transform (PMHT) is limited by its fixed basis functions and single spectral representation, which restricts its ability to capture complex geometric features. This paper proposes a point cloud manifold fractional harmonic transform (PMFHT), which generalizes PMHT by introducing fractional-order parameters and constructs a continuously adjustable intermediate fractional-order spectral domain between the spatial domain and the frequency domain. This fractional-order framework supports more flexible transformation and filtering operations. Experiments show that choosing different transformation orders can enrich the spectral representation of point clouds and achieve excellent results in tasks such as filtering and feature enhancement. Therefore, PMFHT not only expands the theoretical framework of point cloud spectral analysis, but also provides a powerful new tool for manifold geometry processing.
\end{abstract}
\section{Introduction}
\label{sec:intro}

Point clouds are one of the most common data 
representations for 3D scenes, obtained from LiDAR, structured-light scanners, and multi-view stereo reconstruction \cite{rusu20113d}. Unlike meshes, point clouds are unstructured and lack explicit connectivity information, making classical grid-based spectral methods inapplicable \cite{vallet2008spectral,karni2000spectral}. Manifolds are topological spaces that locally exhibit Euclidean properties\cite{li2025novel,hua2023lda}. To enable spectral analysis, point-based manifold harmonic bases (PB-MHB) were introduced \cite{liu2012point}, generalizing the Laplace–Beltrami operator (LBO) to point-sampled manifolds. By solving the discrete LBO eigenvalue problem, the obtained orthogonal basis supports a Fourier-like transform on point clouds, namely the Point Manifold Harmonic Transform (PMHT). This enables applications such as spectral filtering, compression, and feature extraction.

Despite its usefulness, PMHT is inherently limited by its fixed harmonic basis and single spectral representation. The resulting transform only produces one spectrum, which may not fully capture complex geometric features, especially when multi-scale or localized spectral information is required. In signal processing, fractional Fourier transform (FRFT) address a similar limitation by introducing a fractional order parameter\cite{tao2006research,tao2008sampling,wei2016generalized}, providing a continuum of intermediate representations between the time domain and frequency domain. Inspired by this idea, we extend manifold harmonics to a fractional-order setting.

In this paper, we propose the point cloud manifold fractional harmonic transform (PMFHT), a generalization of PMHT. The key idea is to nonlinearly scale the eigenvalues of the LBO using a fractional-order parameter $\alpha$, yielding a continuously tunable family of fractional spectral domains between the spatial and frequency domains. This enables transforms and filtering operations that adapt to different levels of spectral resolution. The mathematical formulation of PMFHT is derived using the fractional power definition based on the eigen-decomposition of the manifold harmonic matrix.
The main contributions of this work are summarized as follows:
\begin{itemize}
  \item We construct a unified framework for fractional harmonic transforms on point cloud manifolds, generalizing the standard PMHT to a fractional-order spectral representation.
  \item We provide a fractional power formulation of PMFHT that is simple, efficient, and theoretically sound.
  \item We demonstrate that PMFHT enables richer spectral representations, leading to superior performance in point cloud denoising, feature enhancement, and shape analysis compared to conventional PMHT.
\end{itemize}

\section{Preliminaries}
\label{sec:format}
\subsection{Point-based manifold harmonic
 transform (PMHT)}
PB-MHB extend the classical manifold harmonics framework to point cloud, enabling spectral analysis without requiring explicit mesh connectivity.
The construction of PB-MHB relies on a symmetrizable discrete Laplace--Beltrami operator (LBO) defined directly on point clouds, ensuring the orthogonality of the resulting bases.


For a compact, boundaryless Riemannian manifold $\mathcal{M}$, the Laplace--Beltrami operator is defined as
\begin{equation}
\Delta_{\mathcal{M}} f = \mathrm{div} \, \mathrm{grad} \, f.
\end{equation}
The eigenvalue problem of the LBO is given by
\begin{equation}
\Delta_{\mathcal{M}} H = -\lambda H,
\end{equation}
where $\lambda \geq 0$ and $H$ are the eigenvalues and eigenfunctions, respectively. The eigenfunctions $\{H_i\}$ form an orthogonal basis with respect to the $L^2$ inner product on $\mathcal{M}$:
\begin{equation}
\langle H_i, H_j \rangle = 0, \quad \text{for } i \neq j.
\end{equation}


Given a point cloud $P = \{p_1,\dots,p_N\}$ sampled from $\mathcal{M}$ satisfying the $\left(\varepsilon, s \varepsilon \right)$-sampling condition, PB-MHB construction begins with a pointwise approximation of $\Delta_{\mathcal{M}} f(p)$. The continuous integral representation of $\Delta_{\mathcal{M}}$ can be written as \cite{liu2012point}
\begin{equation}
\begin{split}
\Delta_{\mathcal{M}} f(p) &= 
\lim_{t \to 0} \frac{1}{4 \pi t^2} \Bigg[
f(p) \int_{\mathcal{M}} e^{-\frac{\|p-y\|^2}{4t}} d\mu(y) \\
&\quad - \int_{\mathcal{M}} e^{-\frac{\|p-y\|^2}{4t}} f(y) d\mu(y)
\Bigg].
\end{split}
\end{equation}

For each point $p \in P$, estimate its local tangent plane $\hat{T}_p$ by performing PCA on its $r$-neighborhood $P_r = P \cap B(p, r)$, where $r = 10 \varepsilon$.
Project neighbors $P_\delta = P \cap B(p,\delta)$, with $\delta \ge 10 \varepsilon$, onto $\hat{T}_p$ and compute their Euclidean Voronoi diagram. The area of the Voronoi cell $\mathrm{vol}( \mathrm{Vor}_{\hat{T}_p}(p) )$ approximates the local area element around $p$ on $\mathcal{M}$.

The discrete approximation of $\Delta_{\mathcal{M}} f(p)$ is then
\begin{equation}
\widehat{\Delta}_t f(p) =
\frac{1}{4 \pi t^2} \sum_{q \in P_\delta}
e^{-\frac{\|q-p\|^2}{4t}}
\bigl(f(q) - f(p)\bigr) \, \mathrm{vol}\!\bigl(\mathrm{Vor}_{\hat{T}_q}(q)\bigr).
\label{eq:discrete-lbo}
\end{equation}
As the sampling becomes denser and $t(\varepsilon) = \varepsilon^{1/2+\alpha}$ with $\alpha>0$, $\widehat{\Delta}_t f(p)$ converges pointwise to $\Delta_{\mathcal{M}} f(p)$.


Collecting \eqref{eq:discrete-lbo} for all $p_i \in P$ yields the matrix form
\begin{equation}
\widehat{\Delta}_t f = L f,
\end{equation}
where $L = B^{-1} Q$ is symmetrizable with
\[
q_{ij} =
\frac{\mathrm{vol}(\mathrm{Vor}_{\hat{T}_{p_i}})\mathrm{vol}(\mathrm{Vor}_{\hat{T}_{p_j}})}{4 \pi t^2}
\exp\!\!\left(-\frac{\|p_i - p_j\|^2}{4t}\right), \quad i \neq j
\]

\[
q_{ii} = - \sum_{j \neq i} q_{ij}
\]

\[
b_{ii} = \mathrm{vol}(\mathrm{Vor}_{\hat{T}_{p_i}})
\]


Solving the generalized eigenproblem
\begin{equation}
QH = \lambda B H
\end{equation}
yields eigenvalues $\{\lambda_i\}$ and orthogonal eigenvectors $\{H_i\}$, the Point-Based Manifold Harmonic Bases (PB-MHB). These bases allow a Fourier-like spectral decomposition of any function $f$ sampled on $P$:
\begin{align}
\tilde{f}_i &= \langle f, H_i \rangle = f^\top B H_i, \\
f &= \sum_i \tilde{f}_i H_i.
\end{align}
This process is known as the PMHT, and is fundamental for spectral filtering, compression, and feature extraction on point clouds.
\subsection{Fractional fourier transform (FRFT)}
\subsubsection{Integral Transform Definition of FRFT}

The FRFT is defined as \cite{alikacsifouglu2024graph}
\[
X_p(u) = F_p[x](u) = \int_{-\infty}^{+\infty} x(t) K_p(t,u) \, dt, \tag{10}
\]
where \(K_p(t,u) = \sqrt{1 - j \cot \alpha} \, e^{\, j \pi ( t^2 \cot \alpha - 2 \sec \alpha + u^2 \cot \alpha )}\),\\
\(\alpha \neq n\pi\); then 
if \(\alpha = 2n\pi\), \(K_p(t,u) = \delta(t-u)\) 
and if \(\alpha = (2n \pm 1)\pi\), \(K_p(t,u) = \delta(t+u)\).
where \(\alpha = p \pi / 2\) indicates the rotation angle of the transformed signal for the FRFT, \(p\) is the transform order, and \(F_p\) denotes the FRFT operator. It is obvious that the FRFT is periodic with period 4. In particular, if \(p = 4n+1\) (\(\alpha = 2n\pi + \pi/2\)), then the FRFT reduces to the conventional Fourier transform. Let \(u = u/\sqrt{2\pi}\) and \(t = t/\sqrt{2\pi}\).

\subsubsection{Fractional Power of FT Definition of FRFT}

The FRFT can also be defined as the $a^{\text{th}}$ fractional power of the ordinary FT \cite{wang2017fractional}:

\[
\mathcal{F}^{a}\psi_{\ell}=(e^{-j\ell\pi/2})^{a}\psi_{\ell}=e^{-j a\ell\pi/2}\psi_{\ell},
\tag{11}
\]
where $\psi_{\ell}$, $\ell=0,1,\ldots$ are Hermite-Gaussian functions, which are the eigenfunctions of the FT operator with eigenvalues $e^{-j\ell\pi/2}$. The Hermite-Gaussian functions are defined as:

\[
\psi_{n}(u)=\frac{2^{1/4}}{\sqrt{2^{n}n!}}H_{n}(\sqrt{2\pi}u)e^{-\pi u^{2}},
\tag{12}
\]
where $H_{n}$ is the $n^{\text{th}}$ order Hermite polynomial. This version defines FRFT by specifying its eigenvalues and eigenfunctions (eigenvectors for the discrete case). This approach ultimately defines the kernel in (10) as the spectral expansion over Hermite-Gaussian basis:

\[
K_{a}(u,u^{\prime})=\sum_{\ell=0}^{\infty}\psi_{\ell}(u)e^{-ja\ell\pi/2}\psi_{\ell}(u^{\prime}),
\tag{13}
\]
and is discretized to define the DFRFT through a transformation matrix $\mathbf{F}^{a}$. Analogous to the spectral expansion of the kernel, where each element of the order-$a$ DFRFT matrix, $\mathbf{F}^{a}$, is defined as follows:

\[
\mathbf{F}_{m,n}^{a}=\sum_{\substack{k=0 \\ k\neq N-1+(N)_{2}}}^{N}\phi_{m}^{(k)}e^{-j\frac{a\pi}{2}k}\phi_{n}^{(k)},
\tag{14}
\]
where $(N)_{2}\equiv N\mod 2$, $\phi_{n}^{(k)}$ is the $n^{\text{th}}$ entry of $k^{\text{th}}$ discrete Hermite-Gaussian vector.

\section{Fractional harmonic transform on point cloud manifolds (PMFHT)}
\label{sec:pagestyle}
\begin{figure}[htb]
  \centering
  \includegraphics[width=0.5\linewidth]{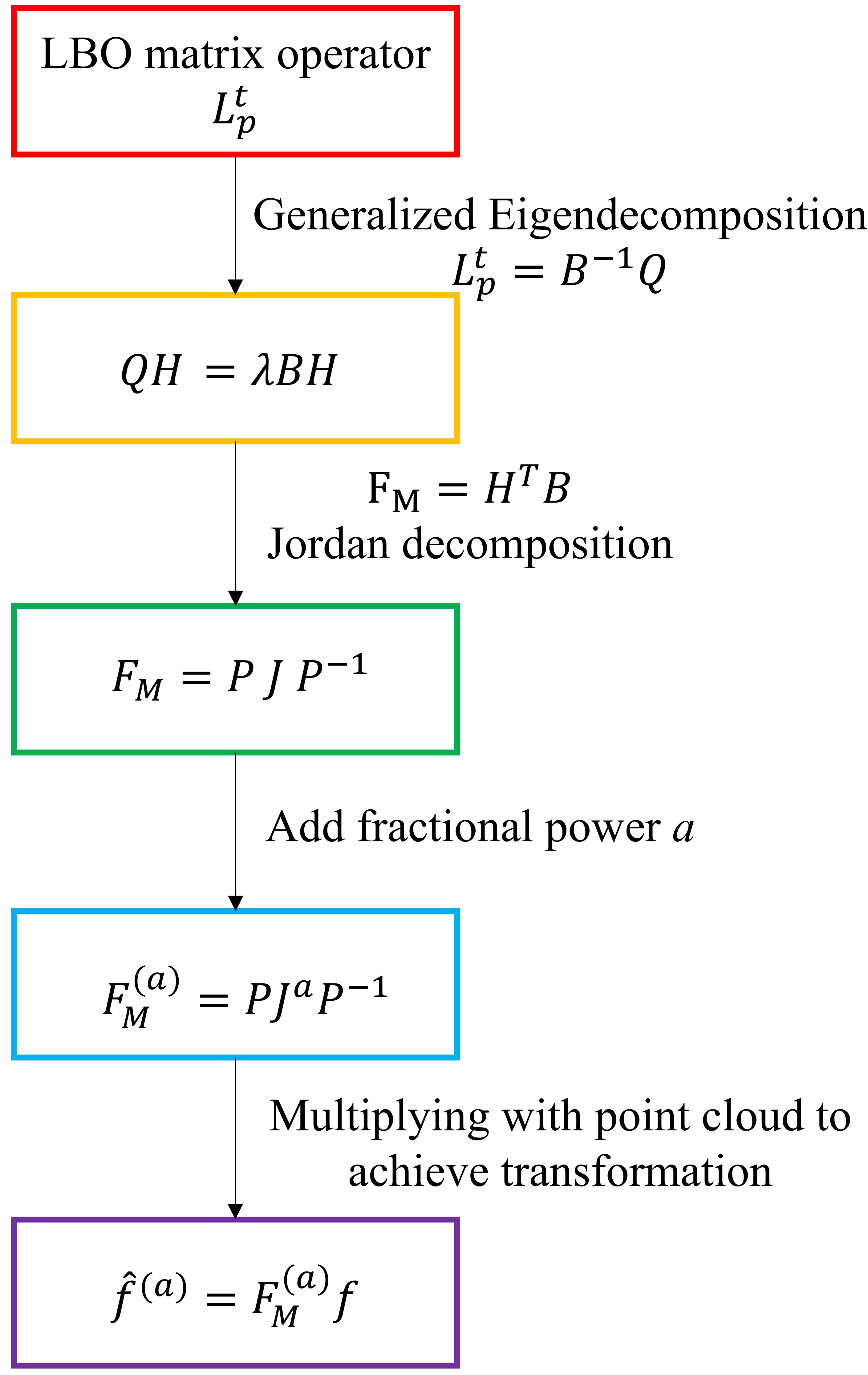}
  \caption{Mathematical principles flowchart}
  \label{fig:liucheng}
\end{figure}

The implementation principle of PMFHT is shown in Fig. ~\ref{fig:liucheng}.
\subsection{Definition}
The complete LBO matrix operator is:
\[
L_P^t = B^{-1}Q, \tag{15}
\]

After obtaining the discrete LBO, solve the generalized eigenvalue problem:
\[
QH = \lambda BH \tag{16}
\]
where the eigenvector set $H$ is called the point cloud manifold harmonic basis.

The forward and inverse point cloud manifold harmonic transforms (PMHT) are defined as:
\[
\hat{f} = H^T Bf, \tag{17}
\]

\[
\quad f = H\hat{f}.\tag{18}
\]
where $H^T BH=I$.

From a matrix operation perspective, the forward transform can be written as:
\[
\hat{f} = F_M f,\tag{19}
\]
where
$F_M := H^T B$
is the point cloud manifold Fourier matrix, and its inverse is $F_M^{-1} = H$.

Let $F_M = P J P^{-1}$, and define the point cloud manifold fractional Fourier matrix of order $a \in \mathbb{R}$ as:
\[
F_M^{(a)} := P J^a P^{-1}.\tag{20}
\]

The point cloud manifold fractional harmonic transform of order $a \in \mathbb{R}$ is then defined as:
\[
\hat{f}^{(a)} = F_M^{(a)} f, \quad f = (F_M^{(a)})^{-1} \hat{f}^{(a)} = F_M^{(-a)} \hat{f}^{(a)}.\tag{21}
\]
\subsection{Properties}

\begin{enumerate}[label=\arabic*)]
  \item \text{Zero rotation} \\
  \begin{equation}
    F_0 = P J_0 P^{-1} = P I P^{-1} = P P^{-1} = I \tag{22}
  \end{equation}
  so the 0th-order PMFHT of a point cloud is the point cloud itself.

  \item \text{Reduction to PMHT when $\alpha = 1$} \\
  \begin{equation}
    F_1 = P J_1 P^{-1} = P J P^{-1} = F_M \tag{23}
  \end{equation}
  so the 1st-order PMFHT of a point cloud is the PMHT of the point cloud.

  \item \text{Index additivity} \\
  As mentioned earlier, the matrix power function is index additive, so we have
\begin{equation}
\begin{split}
F_\alpha F_\beta &= P J_\alpha P^{-1} P J_\beta P^{-1} \\
&= P J_{\alpha + \beta} P^{-1} = F_{\alpha + \beta}.  
\end{split}
\tag{24}
\end{equation}
\end{enumerate}

\section{EXPERIMENT RESULTS}
\label{sec:illust}

\begin{figure}[htb]
  \centering
  \includegraphics[width=0.8\linewidth]{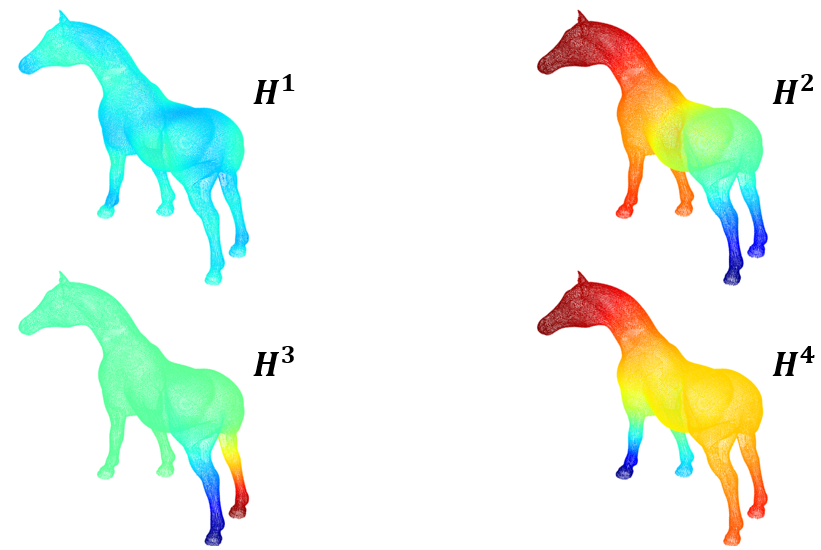}\\
  \includegraphics[width=0.8\linewidth]{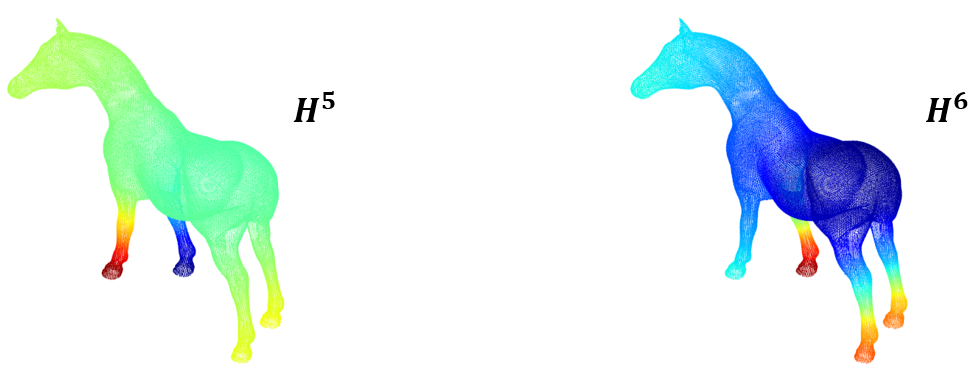}
  \caption{Point cloud manifold harmonic orthogonal bases $H_1-H_6$.}
  \label{fig:combined}
\end{figure}

\begin{figure}[htb]
  \centering
  \includegraphics[width=1\linewidth]{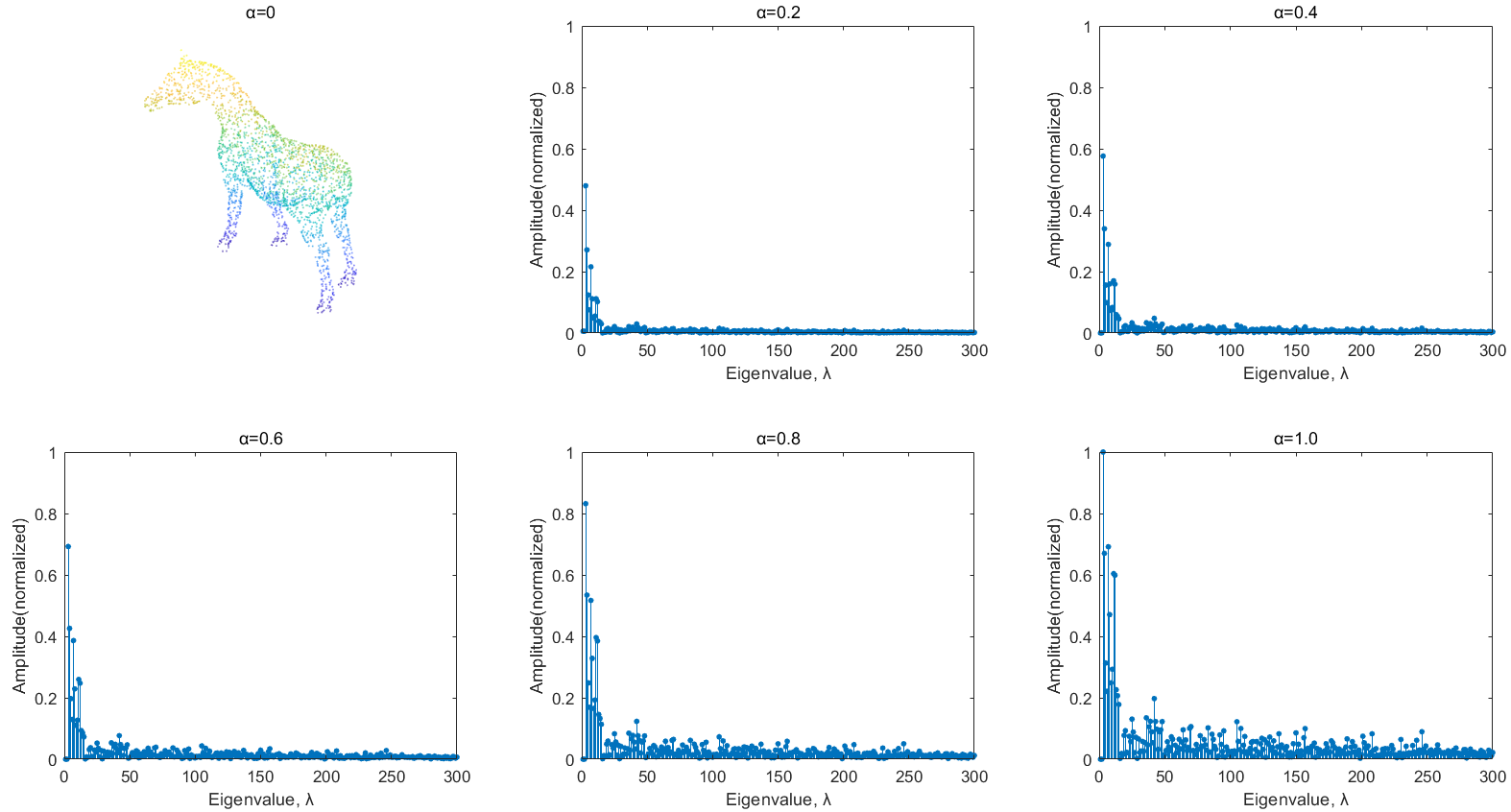}
  \caption{PMFHT spectra of horse point cloud}
  \label{fig:horse}
\end{figure}

\begin{figure}[htb]
  \centering
  \includegraphics[width=1\linewidth]{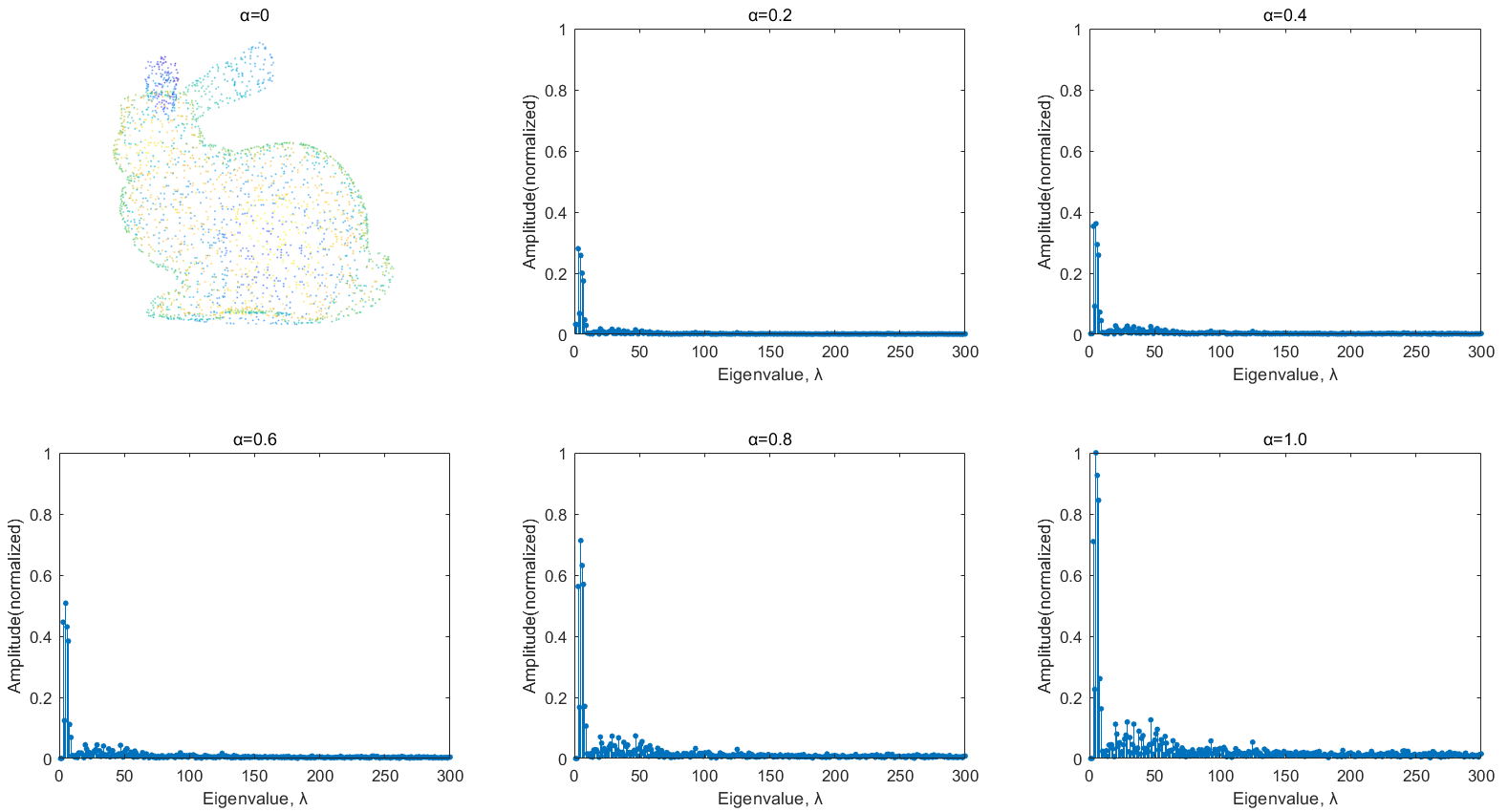}
  \caption{PMFHT spectra of bunny point cloud}
  \label{fig:bunny}
\end{figure}

The point cloud data used in the experiment comes from the Stanford 3D Scan Repository, including horse.ply and bunny.ply. The experimental workflow is as follows: First, the raw point cloud is preprocessed using Python, including loading the point cloud, performing a neighborhood search, and constructing the discrete Laplace-Beltrami operator. Then, a generalized eigenvalue problem is solved to obtain the harmonic basis of the point cloud manifold and perform a fractional-order transform. The results are then imported into MATLAB for visualization and saved as images for comparative analysis. Fig. ~\ref{fig:combined} shows the first six eigenfunctions obtained from the PMHT of the horse model, with the color map corresponding to the eigenfunction value at each vertex. Fig. ~\ref{fig:horse} and ~\ref{fig:bunny} show the results of applying PMHT of different orders to the horse and bunny models, respectively. Experiments show that when the fractional-order parameter $\alpha = 0$, the output of the PMFHT is identical to the original point cloud. However, as $\alpha$ changes, the point cloud signal is gradually projected into different fractional-order spectral domains, resulting in a continuously changing fractional-order energy spectrum, which can flexibly capture multi-scale geometric features.

\begin{figure}[htb]
  \centering
  \begin{subfigure}[b]{0.2\textwidth} 
    \centering
    \includegraphics[width=0.65\linewidth]{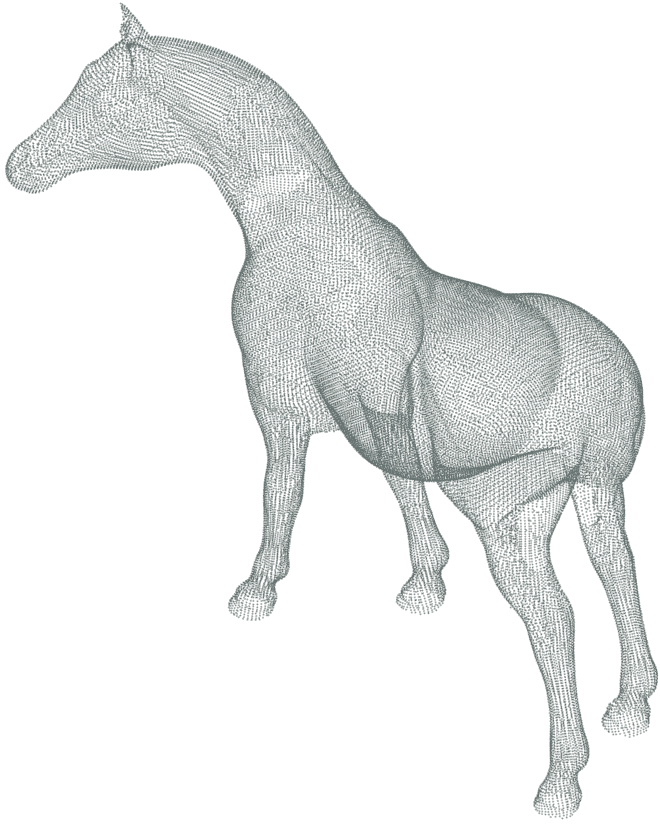}
    \caption{Original point cloud data}
    \label{fig:sub1}
  \end{subfigure}
  
  \vspace{0.15cm} 
  \begin{subfigure}[b]{0.2\textwidth}
    \centering
    \includegraphics[width=0.65\linewidth]{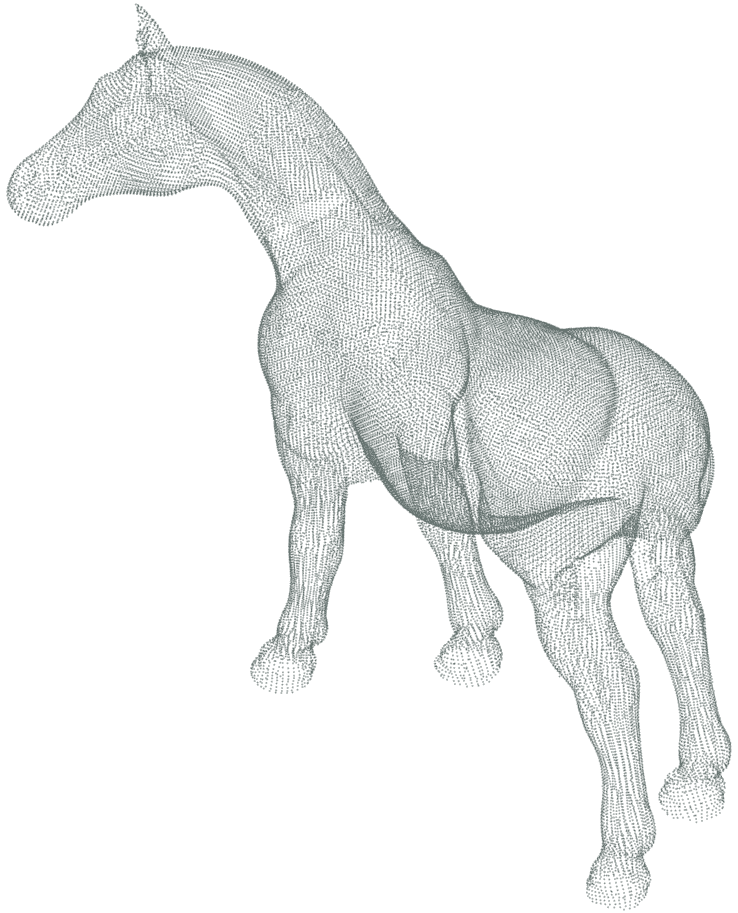}
    \caption{High-pass filtering result}
    \label{fig:sub2}
  \end{subfigure}
  \hfill
  \begin{subfigure}[b]{0.2\textwidth}
    \centering
    \includegraphics[width=0.65\linewidth]{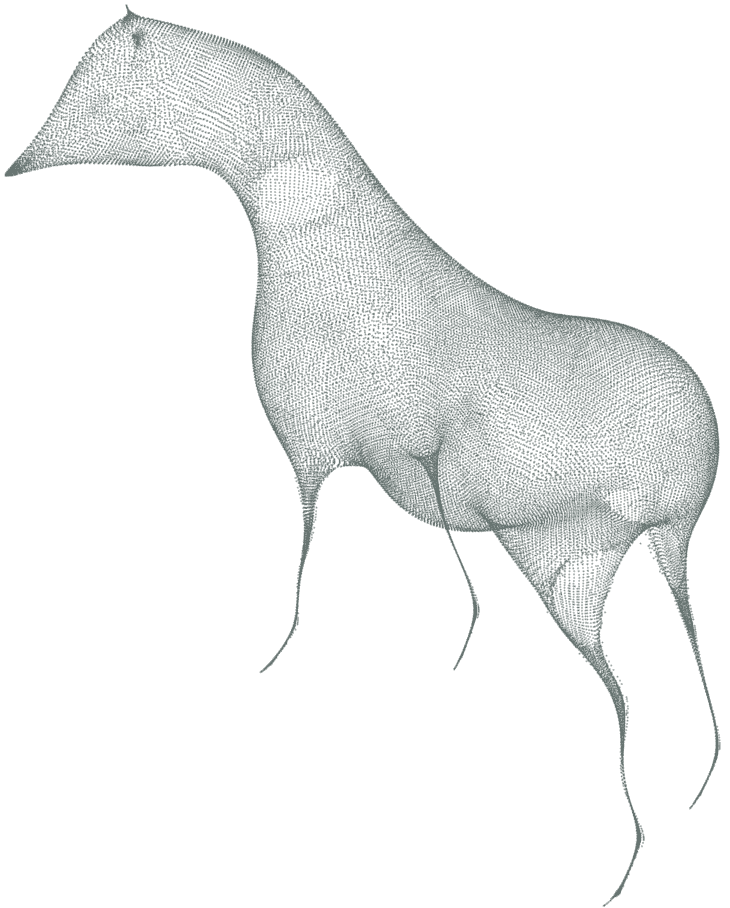}
    \caption{Low-pass filtering result}
    \label{fig:sub3}
  \end{subfigure}

  \caption{Point cloud manifold fractional harmonic filtering results.}
  \label{fig:overall}
\end{figure}

Furthermore, Fig. ~\ref{fig:overall} demonstrates the effect of spectral filtering using PMFHT. By selecting the appropriate fractional order and filter type, high-pass filtering and low-pass filtering can be implemented, respectively. The former enhances high-frequency details in the point cloud, such as edges and sharp geometric features, while the latter smooths noise and preserves low-frequency shape information.

\section{CONCLUSIONS}
\label{sec:cons}
This paper presents PMFHT, a generalized spectral analysis framework that extends PMHT with fractional-order parameters. PMFHT offers continuous, adjustable fractional-order spectral representations, bridging spatial and frequency domains. This allows more adaptive decomposition of point cloud data. Experiments on real point clouds show that PMFHT can enhance or smooth geometric details via fractional-order filtering, providing richer information than standard PMHT. By combining fractional-order signal processing with point cloud spectral geometry, it broadens the theoretical basis and offers a powerful tool for point cloud processing.

\section{ACKNOWLEDGEMENT}
\label{ackn}
This work was supported by grants from the National Natural Science Foundation of China [No. 62171041].



\clearpage
\bibliographystyle{IEEEbib}
\bibliography{strings,refs}

\end{document}